\documentclass{article}

\usepackage[a4paper, total={6in, 9in}]{geometry}

\usepackage{amsfonts}
\usepackage{amsmath}
\usepackage{amssymb}
\usepackage{amsthm}
\usepackage{graphics}
\usepackage{graphicx}
\usepackage{authblk}
\usepackage[colorlinks,breaklinks,backref]{hyperref}
\usepackage{hyperref}
\usepackage{backref}

\usepackage[normalem]{ulem}
\usepackage{color}
\definecolor{gold}{rgb}{0.55,0.22,0.0}

\def\proof{\par\noindent{\bf Proof\ \ }}

\def\qed{\hfill $\square$ \\}

\def\Aut{{\rm Aut}}

\newtheorem{theorem}{Theorem}[section]
\newtheorem{lemma}[theorem]{Lemma}
\newtheorem{corollary}[theorem]{Corollary}
\newtheorem{proposition}[theorem]{Proposition}

\begin{document}

\baselineskip15pt

\title{Distinguishing graphs of maximum valence 3}
\author[1]{Svenja H\"uning}
\author[2]{Wilfried Imrich}
\author[1]{Judith Kloas}
\author[1]{Hannah Schreiber}
\author[3]{Thomas W. Tucker\thanks{Supported by the Austrian Science Fund FWF W1230 and P24028 and by Award 317689 from the Simons Foundation.}}
\affil[1]{
	\footnotesize Graz University of Technology, 8010 Graz, Austria \protect \\ 
	e-mail: huening@tugraz.at, kloas@math.tugraz.at, hschreiber@tugraz.at
}
\affil[2]{
	\footnotesize Montanuniversit\"at Leoben, 8700 Leoben, Austria \protect \\ 
	e-mail: imrich@unileoben.ac.at
}
\affil[3]{
	\footnotesize Department of Mathematics, Colgate University, Hamilton, NY, USA \protect \\ 
	e-mail: ttucker@colgate.edu
}
\setcounter{Maxaffil}{0}
\renewcommand\Affilfont{\itshape\small}
\maketitle

\begin{abstract}
The distinguishing number $D(G)$ of a graph $G$ is the smallest number of colors that is needed to color $G$
such that the only color preserving automorphism is the identity. We give a complete classification  for all connected graphs $G$ of maximum valence $\triangle(G)=3$ and distinguishing number $D(G) = 3$. As one of the consequences we get that all infinite connected graphs with $\triangle(G)=3$ are 2-distinguishable.
\end{abstract}

\section{Introduction}\label{sec:intro}

The distinguishing number of a group $A$ acting faithfully on a set $\Omega$ is the least number of colors needed to color the elements of $\Omega$ such that the only color-preserving element of $A$ is the one that fixes all elements of $\Omega$.
If $A$ is the automorphism group of a graph $G$, then the distinguishing number $D(G)$ of $G$ is the distinguishing number of the action of $A$ on the vertex set of $G$. Since its introduction by Albertson and Collins \cite{AC96} more than 20 years ago, there has developed an extensive literature on this topic.

Actually Babai \cite{B77} showed already 1977 that a tree has a distinguishing coloring with two colors if all vertices have the same valence $\alpha\geq2$, where $\alpha$ can be an arbitrary finite or infinite cardinal\footnote{If $\alpha$ is smaller than the first uncountable inaccessible cardinal, then there also exists a coloring with a finite number of colors that is only preserved by the identity endomorphism.}, but the subject lay dormant until the seminal paper of Albertson and Collins \cite{AC96}.

The concept also has had an independent separate history in the theory of permutation groups \cite{CNS84}, unknown to graph theorists until recently \cite{BC11}.

The first motivation for this paper is a bound by Collins and Trenk \cite{CT06} and, independently, Klav\v{z}ar, Wong and Zhu \cite{KWZ06}. They proved that for any finite graph $G$ of maximum valence $\Delta(G) = d$, $D(G) \leq d+1$ with equality only if $G$ is the complete graph $K_{d+1}$, the complete bipartite graph $K_{d,d}$, or the $C_5$. For infinite graphs the bound is the supremum of the valences, see Imrich, Klav\v{z}ar and Trofimov \cite{IKT07}. Hence, for infinite graphs $D(G)\leq d$ if $G$ has bounded valence $d$. If one wishes to improve this bound, it is reasonable to begin with $d = 3$.

The second, equally important motivation, is the Infinite Motion Conjecture of Tucker \cite{T11}, who conjectured that each connected, locally finite infinite graph is 2-distinguishable if every automorphism  that is not the identity moves infinitely many vertices. The conjecture is still open, although it has been shown to be true for many classes of graphs \cite{CT11,ISTW15,STW12}, in particular for graphs of subexponential growth \cite{Lehner2}, and thus for all graphs of polynomial growth. For a long time it was not clear whether it holds for graphs of maximal valence 3, and whether infinite motion was really needed. This was first solved under the additional condition of vertex transitivity \cite{ILTW-xx}. It turns out that all finite or infinite connected, vertex transitive graphs are 2-distinguishable unless they are one of four exceptional graphs.

Here the result is extended to a complete classification of all finite or infinite  connected graphs of maximal valence 3 that are not 2-distinguishable.

We begin with a general observation about graphs of bounded valence.

For any graph $G$ with $\Delta(G) = d$ and $D(G) = d - 1$, one can subdivide an edge with a vertex $v$ and add an edge between $v$ and a vertex of a disjoint copy of $K_d$ to get a graph $G'$ with $\Delta(G') = d$ and $D(G') = d - 1$ (if $G$ is $d$-valent, then $G'$ can be as well, simply by attaching $d - 2$ copies of $K_d$). Thus, the only cases where one might expect a classification of graphs with a given distinguishing number are $D(G) = d$.

There are infinitely many graphs with $\Delta(G) = d$ and $D(G) = d$. Let $T(n,d)$ be the tree where all vertices have valence 1 or $d$ and every vertex of valence 1 has the same distance $n$ from a root vertex $v$. Clearly, $D(T(n,d)) \geq d - 1$ and $D(T(n,d)) = d - 1$ if and only if $D(T(n+1,d)) = d - 1$.  But $T(1,d) = K_{1,d}$ so $D(T(1,d)) = d$ and hence $D(T(n,d)) = d$.

From now on we assume that the maximum valence is 3, unless otherwise stated.

We call a vertex of valence 1 a \emph{leaf}.

If $u$ and $v$ have a common neighbor, we say they are \emph{siblings} or a \emph{sibling pair}.
Vertex $v$ is an \emph{only child} of a vertex $u$ if $v$ is the only neighbor of $u$ with $valence(v) = 1$. For $d = 3$, we abbreviate $T_n = T(n,d)$.

We give some variations of the trees $T_n$ which also have distinguishing number 3. The most obvious one is simply to join each sibling pair of leaves by an edge. Denote this graph $S_n$. We do this only for $n > 1$ since $T_1$ has three sibling pairs and adding all such edges gives $K_4$. Note that we can also think of $S_n$ as obtained from $T_n$ by attaching a triangle to each leaf of $T_{n-1}$.

The other three variation are obtained by replacing the edges between sibling pairs in $S_n$ by three other ``gadgets". In each case, the sibling pair vertices are labeled $u,v$.

\begin{itemize}
	\item \emph{Gadget of type 1}: A $4$-cycle $uxvy$ with $x,y$ valence 2.
	\begin{center}
		\includegraphics[scale=0.8]{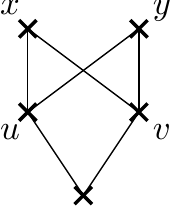}
	\end{center}
	
	\item \emph{Gadget of type 2}: The same as the type 1 gadget but with an edge $xy$.
	\begin{center}
		\includegraphics[scale=0.8]{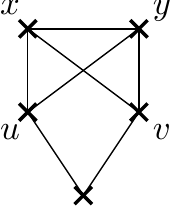}
	\end{center}
	
	\item \emph{Gadget of type 3}: A hexagon $uxzvyw$ with edges $xy$, $zw$ (this can be viewed as $K_{2,2}$ with $u$ joined to one part and $v$ to the other part).
	\begin{center}
		\includegraphics[scale=0.8]{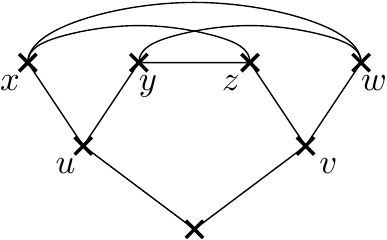}
	\end{center}
\end{itemize}

\begin{figure}[h]
	\centering
	\includegraphics[scale=0.8]{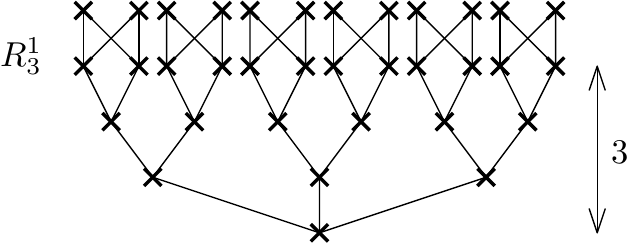}
	\includegraphics[scale=0.8]{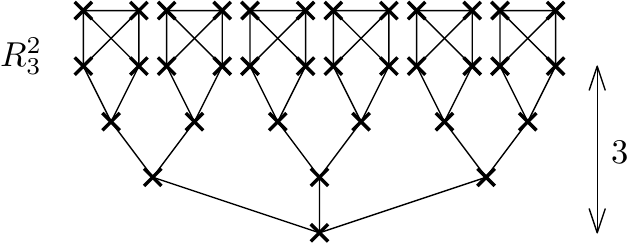}
	\includegraphics[scale=0.8]{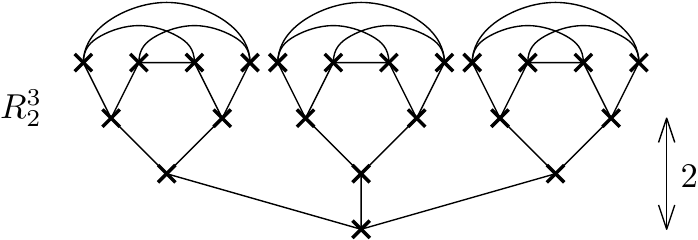}
	
	\caption{$R^1_3$ (top left), $R^2_3$ (top right) and $R^3_2$ (bottom).} \label{fig:tools}
\end{figure}

Now, we define three graphs, $R^1_n$, $R^2_n$ and $R^3_n$, by adding the respective gadgets between each sibling pair of leaves in the tree $T_n$. See Figure~\ref{fig:tools} for examples.

\bigskip

Since $Aut(S_n)$ acts on its vertices the same way as $Aut(T_n)$, we have $D(S_n) = 3$. If there existed a distinguishing 2-coloring for one of the graphs $R_n^1$, $R_n^2$ or $R_n^3$, then this would induce a distinguishing 2-coloring of the associated $T_n$. Therefore $D(G) = 3$ for $G = R_n^1, R_n^2, R_n^3$.

Our classification of graphs with $\Delta(G) = 3$ and $D(G) = 3$ is the following:

\begin{theorem} [Main Theorem] \label{thm:main}
	Let $G$ be a finite or infinite connected graph with $\Delta(G) = 3$. Then $D(G) = 3$ if and only if $G$ is either 
	$K_{1,3}$, $K_{2,3}$, the cube $Q$, the Petersen graph $P$, or a member of one of the five families $T_n$, $S_n$, $R_n^1$, $R_n^2$, $R_n^3$ for $n > 1$.
\end{theorem}

We note for the four exceptions, clearly $D(K_{1,3}) = D(K_{2,3}) = 3$, and it is an exercise to verify that $D(Q) \neq 2$ (or see \cite{T11}). It is slightly more work to show $D(P)\neq 2$; we will sketch a proof in Section 2.

The proof that these are the only graphs $G$ with maximum valence $d = 3$ and $D(G) = 3$ occupies most of the rest of this paper. In Section 2, we first give some Corollaries that may shed some light on the general problem when $d > 3$. In Section 3 we introduce a 2-coloring which either is distinguishing or leads to restrictions on the local structure of $G$. This coloring is then used throughout the rest of the paper. In Section 4, we show that if $G$ has any leaves, then $D(G) = 2$ unless $G = T_n$. The goal of Section 5 is to reduce to the case of edge transitive graphs by analyzing how the stabilizer of a vertex $v$ acts on the neighbors of $v$. In Section 6, we show that any cubic graph $G$ with girth at least 6 has $D(G) = 2$; this proof does not use edge transitivity. This completes the proof of the Main Theorem, since the five edge transitive cubic graphs of girth less than 6 are easily analyzed. In Section 7 we pose a variety of questions.

\section{Corollaries}

We give some corollaries of the Main Theorem, mostly just observations about our list of graphs with $D(G) = 3$. Each gives some insight into the relationship between distinguishing number and graph structure. Each suggests ways one might generalize the case of maximum valence 3 to graphs of higher valence.

\begin{corollary} \label{cor:inf2dist}
	Every infinite connected graph of maximal valence 3 is 2-distinguishable.
\end{corollary}

\begin {corollary} \label{cor:vt}  
	Every  connected, vertex transitive graph of maximal valence 3 is 2-distinguishable, except for $K_{3,3}$, $K_4$, $Q$ and $P$.
\end{corollary}

\begin {corollary} \label{cor:et}  
	Every connected, edge transitive graph of maximal valence 3 that is not vertex transitive is 2-distinguishable, except for $K_{1,3}$ and $K_{2,3}$.
\end{corollary}

For more direct proofs of Corollary \ref{cor:vt} and \ref{cor:et} see \cite{ILTW-xx}.

\begin{corollary} \label{cor:conn}
	Every 2-connected graph of maximal valence 3 is 2-distinguishable, except for $K_{2,3}$, $K_{3,3}$, $K_4$, $Q$, and $P$.
\end{corollary}

The length of the shortest cycle in $G$ is its \emph{girth}. The following result is in fact one of the steps in the proof of the Main Theorem.

\begin{corollary} \label{cor:girth}  
	If $G$ has girth at least 6, then $D(G) = 2$.
\end{corollary}

The \emph{motion} of a group $A$ acting on a set $\Omega$ is the smallest integer $m$ such that some element of $A$ moves exactly $m$ points. The motion of a graph $G$, which we denote $m(G)$, is the motion of $Aut(G)$ acting on the vertex set. The Motion Lemma  \cite{CT11, RS98} states that if $m > 2 \log_2(|A|)$, then the action has distinguishing number 2; the proof is elementary and short. Thus large enough motion gives 2-distinguishability.  For graphs of maximum valence 3, large enough means $3$ or more, except for $Q$ and $P$, since it is easily checked that all other $G$ in our Main Theorem have motion 2.

\begin{corollary} \label{cor:motion}
	If $m(G) > 2$, then $D(G) = 2$ with the exception of $Q$ and $P$.
\end{corollary}

In fact, when $D(G) = 3$ and $G$ is not $Q$ or $P$, we can isolate an automorphism of motion 2 using a 2-coloring of $G$. We say a coloring \emph{fixes} a set of vertices if any color-preserving automorphism fixes all vertices in that set.

\begin{corollary} \label{cor:wi}
	 If $D(G) = 3$ and $G$ is not $Q$ or $P$, then $G$ admits a 2-coloring that fixes all vertices except two siblings. 
\end{corollary}

\proof
	Clearly, such a 2-coloring exists for $K_{1,3}$ and $K_{2,3}$.
	All other graphs that satisfy the assumptions of the lemma have a root vertex, say $v_0$, corresponding to the root of $T_n$. The $T_n$ are the only graphs in the class with leaves and they all come in sibling pairs. We first construct the desired coloring for $T_n$.

	We begin with a 2-coloring of $T_1$: we color $v_0$ black, two of its neighbors white and one black. Clearly this coloring fixes all vertices except for one pair of interchangeable siblings. To color $T_n$ we color its subgraph $T_1$ as before, and continue inductively by assigning different colors to any two vertices of distance $k > 1$ from $v_0$ if they have a common neighbor of distance $k-1$ from $v_0$. When $k = n$ we make an exception for a single sibling pair of vertices whose shortest paths to $v_0$ contain a white neighbor of $v_0$.
	Both vertices in that pair are colored white. It is easy to see that this coloring fixes all vertices not in this pair.

	For $S_n$, $R_n^1$, $R_n^2$ we proceed analogously, and let the vertices $x,y$ in the gadgets play the role of the sibling pairs in $T_n$. For $R_n^3$ we assign different colors to all gadget vertices $z,w$, but treat the pairs $x,y$ as before. Again, our coloring fixes all vertices except the ones in the white $(x,y)$-pair.
\qed
 
It is easily verified that $m(Q) = 4$.  
 
\begin{proposition}\label{lem:P}
	For the Petersen graph, $m(P) = 6$ and $D(P) = 3$.
\end{proposition}

\proof 
	We follow a remark of Lehner\footnote{Private communication.}. He observed that $P$ is the complement of the line graph of $K_5$. Thus any  edge coloring of $K_5$ is a vertex coloring of the complement $\overline{P}$ of $P$, and thus also of $P$.
	Because every 2-edge coloring of $K_5$ corresponds to a subgraph of $K_5$ and its complement, every subgraphs of $K_5$, together with its complement, yields correspond a 2-colorings of $P$.
	Furthermore, given a subgraph, say $H$, of $K_5$ the group induced by $\Aut(H)$ on $E(H)$ is the same as the group that preserves the coloring of the vertices of $P$ induced by $H$.

	As  the smallest graph with trivial automorphism group has a least six vertices, every subgraph $H$ of $K_5$ has a non-trivial automorphism $a$. We consider the cycle decomposition of $a$. If it has a five-cycle or a four-cycle, then $\Aut(H)$ moves all ten pairs of distinct vertices of $K_5$. If it is a three-cycle it moves at least nine pairs. If $a$ has only 2-cycles or fixed points it moves at least six pairs (but can fix four).

	For a 3-distinguishing coloring of $P$ we let $H$ be a path of length 4 in $K_5$ and choose an end-edge $e$ of $H$. We color $e$ red, the other edges of $H$ black and the edges of the complement of $H$ white. This yields a 3-coloring of $P$. As the group of $H$ that preserves its edge-coloring is trivial, this is a distinguishing 3-coloring of $P$.
\qed

\section{Canonical 2-colorings rooted at a subgraph}

Let $G$ be a cubic or subcubic graph  and  $K$ be a vertex-induced, connected subgraph
with at least one \emph{internal} vertex, that is, a vertex  all of whose neighbors are in $K$.
Define  $S_n(K)$ as the set of vertices of distance $n$ from $K$; one might call it the sphere of radius $n$ about $K$. Thus $S_0(K) = K$ and $S_1(K)$ is the set of vertices not in $K$ but adjacent to some vertex in $K$.
Let $B_n(K)$ denote all vertices of distance at most $n$ from $K$ (the ball of radius $n$ about $K$). For a vertex $v$ in $S_n(K)$, we call its neighbors in $S_{n+1}(K)$, $S_n(K)$, $S_{n-1}(K)$, respectively, its {\em up, cross, down} neighbors. Notice that all vertices which are not in $K$ have at least one down neighbor and that not every vertex has to have an up neighbor.

The idea of constructing a 2-distinguishing coloring of the vertices of $G$ is to color all vertices of $K$ black and then to extend the coloring inductively from one $S_n(K)$ to the next. Our objective is to obtain a 2-coloring of $G$ such that the only color-preserving automorphism fixing the vertices of $K$ is the identity. Thus at stage $n$ we have this:

{\bf Goal} \emph{Assume $B_n(K)$, $n > 0$, has been 2-colored so that any automorphism of $G$ fixing the vertices of $K$ and preserving the coloring of $B_n(K)$ is the identity on $B_n(K)$. Then extend this to a 2-coloring of $S_{n+1}(K)$ that has the same property on $B_{n+1}(K)$.}

The plan for extending the coloring to $S_{n+1}(K)$ is simple enough: if a vertex $v$ in $S_n(K)$ has a single up neighbor, color it white, and if it has two up neighbors that can be switched by an automorphism of $G$ that fixes $S_n(K)$, color one white and one black. The problem is that the up vertices of $v$ may have already been colored when we colored the up neighbors of a different vertex. In the following three paragraphs we will make this procedure more precise.

Assume we have colored the graph up to the sphere $S(n,v)$. Let $V$ be the set of vertices of $S(n, v)$ and $U$ be the set of vertices of $S(n+1,v)$. Moreover, let $H$ be the subgraph of $G$ determined by all edges between vertices of $V$ and $U$. In fact, these are the up-edges from $V$. Note that for the moment we do not care about possible cross-edges in $U$ or $V$. We are interested in coloring the vertices of $U$ such that any automorphism of $G$ which fixes the vertices of $V$ and preserves the coloring of $U$ also fixes the vertices of $U$.

Suppose that $U$ has a vertex $x$ of valence 3 in $H$. If there is another vertex $y \in U$ adjacent to the same vertices in $V$ as $x$, color $x$ black and $y$ white. If there is no such vertex $y$, color $x$ white.

\begin{figure}[h]
	\centering
	
	\includegraphics[scale=0.9]{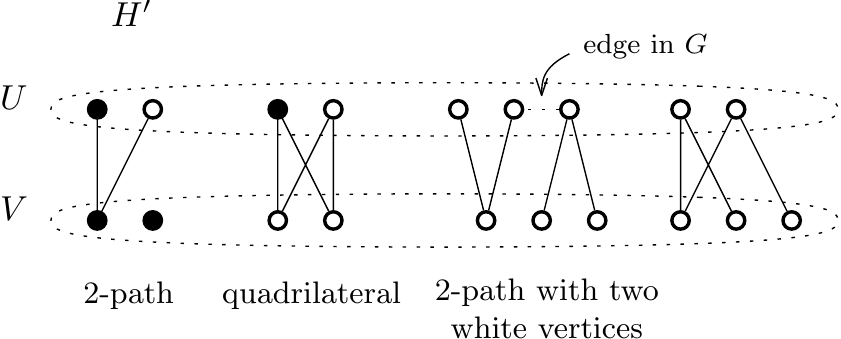}
	
	\caption{Possible configurations in $H'$ used in the description of the canonical coloring.} \label{fig:canonical_coloring}
\end{figure}

Now, consider the subgraph $H'$ of $H$ obtained by removing all valence 3 vertices of $U$ in $H$, i.e. we remove all colored vertices of $U$. The remaining subgraph $H'$ contains only vertices which are of valence 1 or 2, so it is a union of paths and/or cycles such that the vertices in each component alternate between $U$ and $V$. See Figure~\ref{fig:canonical_coloring} for some examples. By assumption the vertices of $V$ are fixed. Therefore, any component of $H'$ is fixed except for two configurations. One is a 2-path, which consists of three vertices such that the middle vertex is in $V$, and the other one is quadrilateral, see Figure~\ref{fig:canonical_coloring}. For all other possible components of $H'$ there is no automorphism of $G$ fixing $V$ but acting non trivial on the considered component of $U$.
Thus, we color all vertices which are neither in a 2-path nor in a quadrilateral white.
For the remaining pairs $x,y \in V$ in a 2-path or quadrilateral, choose a coloring using black the fewest times such that any automorphism of $G$ fixing $U$ also fixes these remaining pairs. Note that we know there are such colorings because we could simply color each such pair black-white.

We call this a {\em canonical 2-coloring} of $G$ {\em rooted at $K$}.

In what follows, we use subscripts to denote which sphere a vertex is in: so $u_n, v_n, x_n \dots$ are vertices in $S_n(K)$. We also extend the notion of siblings as we have defined it in Section \ref{sec:intro} by calling two (distinct) vertices $u_n, v_n$  \emph{siblings} if they have a common down neighbor.

\medskip

We make the following observations about the resulting coloring:

\begin{proposition}[White Up] \label{prop:white_up}
	If $v_n$ has an up neighbor, it has a white up neighbor.
\end{proposition}

\begin{proposition}[Black-white Siblings] \label{prop:black_white_sibling}
	If $v_n$, $n > 0$, is black, it has a white sibling $u_n$ and there is an automorphism of $G$ interchanging $v_n$ and $u_n$, but fixing all other vertices of $B_n(K)$.
\end{proposition}

\proof
	Suppose there was no such automorphism. Then we could color $v_n$ white, contradicting the minimality in the use of black.
\qed

\begin{proposition}[Black Cross] \label{prop:black_cross}
	Suppose $n > 0$. If $u_n$ and $v_n$ are both black and adjacent, then there is a quadrilateral $u_n v_n x_n y_n$, where $x_n$ is a white sibling of $u_n$ and $y_n$ a white sibling of $v_n$.
\end{proposition}

\proof
	Since $v_n$ is black, it has a white sibling $y_n$ with an automorphism interchanging $v_n$ and $y_n$ but fixing $u_n$, forcing an edge $u_ny_n$. Similarly, since $u_n$ is also black there is a white vertex $x_n$ and edge $v_nx_n$. Since the interchange of $v_n$ and $y_n$ also leaves $x_n$ fixed, which is adjacent to $v_n$, we have $y_n$ adjacent to $x_n$.
\qed

\begin{proposition}[All Black] \label{prop:all_black}
	There is no vertex $v_n$ for $n > 1$ of valence 2 or 3 that is black with all neighbors black.
\end{proposition}

\begin{figure}[h]
		\centering
		
		\includegraphics[scale=0.9]{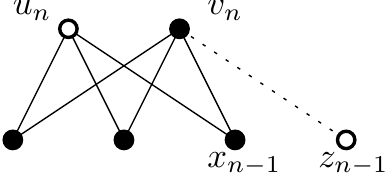}
		
		\caption{The vertex $v_n$ has already three neighbors. Therefore there can not be an edge between $v_n$ and $z_{n-1}$.} \label{fig:all_black}
	\end{figure}

\proof
	By Propositions White Up (\ref{prop:white_up}) and Black Cross (\ref{prop:black_cross}), all neighbors of $v_n$ are down neighbors. Since $v_n$ is black, it has a white sibling $u_n$ (with the same valence) and an automorphism interchanging $u_n$ and $v_n$ fixing $S_{n-1}(K)$. Since all the neighbors of $v_n$ are down neighbors, $u_n$ has the same neighbors as $v_n$. Let $x_{n-1}$ be one of the common neighbors. Since it is black, there is an automorphism $\phi$ interchanging $x_{n-1}$ with the white vertex $z_{n-1}$ and leaving fixed all other vertices of $S_{n-1}(K)$. Since the other down neighbors of $u_n, v_n$ are all black, $\phi$ fixes all of them. Since each of these other possible down neighbors has $u_n, v_n$ as up neighbors, $\phi$ either fixes $u_n, v_n$ or interchanges them. This forces another edge from either $u_n$ or $v_n$ to the white vertex $z_{n-1}$, a contradiction, see Figure \ref{fig:all_black}.
	

\qed

Call a vertex of a canonical coloring {\em kiwi} if it is black surrounded by black. Proposition All Black (\ref{prop:all_black}) says the only kiwi vertices are in $K \cup S_1(K)$.

\begin{proposition} [Internal] \label{prop:internal}
	The only non-identity color-preserving automorphisms of a canonical coloring rooted at the subgraph $K$ are those taking an internal vertex of $K$ to either an internal vertex of $K$ or a kiwi vertex of $S_1(K)$.
\end{proposition}

Note that the neighbors of a kiwi vertex in $S_1(K)$ must be in $K$ by Propositions White Up (\ref{prop:white_up}) and Black Cross (\ref{prop:black_cross}).

\section{Leaves}

We first show that the only graphs $G$ having a leaf with $D(G)=3$ are the trees $T_n$.

\begin{theorem} \label{thm:leaves}
	If $G$ has a leaf, then $D(G) = 2$ or $G = T_n$ for some $n$.
\end{theorem}

\proof
	The smallest two subcubic graphs with a leaf are the $T_1$ and a triangle where one of the three vertices has a further neighbor. For these two graphs the theorem holds.
	Consider now a graph $G$ with more than 4 vertices.
	Suppose some vertex $v$ in $G$ has valence 1 and is an only child of $u$. The canonical coloring where $K$ is the edge $uv$ breaks all automorphisms in $\Aut_v(G)$ since any automorphism fixing $v$ fixes $u$ as well. Since all only children other than $v$ are colored white in a canonical coloring, there is no automorphism moving $v$, so the canonical coloring is distinguishing.

	Now suppose that all leaves of $G$ come in sibling pairs and assume first that $G$ is finite.
	Prune all such sibling pairs. Since $\Delta(G) = 3$, there is at least one leaf in the new graph.
	We have two cases. The new graph has again sibling pairs. Then prune again all sibling pairs and continue like this inductively until you get $T_1$ or a graph with an only child. If we ended with $T_1$ we know that $G = T_n$ for some $n$. 
	Else we get at some point a graph $G'$ with an only child.
	Give $G'$ a distinguishing 2-coloring and color all sibling pairs, one black and one white. Any automorphism $\phi$ of $G$ takes $G'$ to $G'$, so if $\phi$ is color-preserving, $\phi|_{G'}$ is the identity. But then $\phi$ is the identity on $G$ since all successive removed sibling pairs are colored black-white.

	Now assume that $G$ is infinite. As before, if $G$ has an only child leaf, then $D(G) = 2$. Suppose instead all leaves come in sibling pairs.  We want to prune all such pairs to form a graph $G'$. First assume that $G$ has a cycle $C$. Now we can  induct on the distance $s$ from $C$ to the closest leaf. If $s = 1$, the closest leaf is an only child, since the parent has valence 2 on the cycle $C$, so $D(G) = 2$. Assume $D(G) = 2$ for all infinite graphs with a leaf of distance $s = n$ from a cycle. Then for $s = n + 1$, either $G$ has an only child leaf, or pruning all sibling pairs gives a graph $G'$ with $s = n$. In either case, $D(G) = 2$.

	The same argument applies if we replace cycles by 2-way infinite paths. Thus there remain only trees without 2-way infinite paths, but
	1-way infinite paths. In this case there must be a maximal infinite path whose origin is an only child leaf. To see this, let $R$ be a 1-way infinite path. The set of all 1-way infinite paths that contain $R$ is partially ordered by inclusion and every totally ordered subset has a maximal element, namely its union. By Zorn's Lemma there must be a maximal element. It cannot be  a 2-way infinite path. Hence the maximal element is a 1-way infinite path. Because it is maximal its origin must be an only child leaf.
\qed

\begin{lemma} \label{lem:sibling}
	If $G$ has adjacent sibling vertices of valence 2, then  $D(G) = 2$ or $G = S_n$. If $G$ has a gadget of type 1, 2, or 3, respectively, then $D(G) = 2$ or $G = R_n^1, R_n^2, R_n^3$, respectively.
\end{lemma}

\proof
	Suppose that $G$ has adjacent siblings of valence 2. Thus we know that $G$ cannot be $T_n$. If $G$ has a vertex of valence 1, then $D(G) = 2$ by Theorem~\ref{thm:leaves}. Let $G'$ be the graph obtained by removing all edges between sibling pairs of valence 2. We note that $\Aut(G)$ is a subgroup of $\Aut(G')$. Thus if $D(G') = 2$, we have $D(G) = 2$. Suppose instead that $D(G') = 3$.

	Since $G'$ has vertices of valence 1, $G' = T_n$ by Theorem~\ref{thm:leaves}. Since $G$ has no vertex of valence 1, every vertex of valence 1 in $G'$ comes from the removal of an edge between a sibling pair in $G$, making $G = S_n$.

	The proof for gadgets of type 1, 2, 3 is the same, where $G'$ is obtained by removing all gadgets of one type, creating vertices of valence 1.  Any distinguishing 2-coloring of $G'$ extends to one of $G$ by coloring $u_n, v_n$ and $x_n, y_n$ black white, and for gadgets 3, $z_n, w_n$ black-white. Otherwise, $G = R_n^1, R_n^2, R_n^3$.
\qed

\section{Vertex types}

\begin{figure}[h]
	\centering

	\includegraphics[scale=0.9]{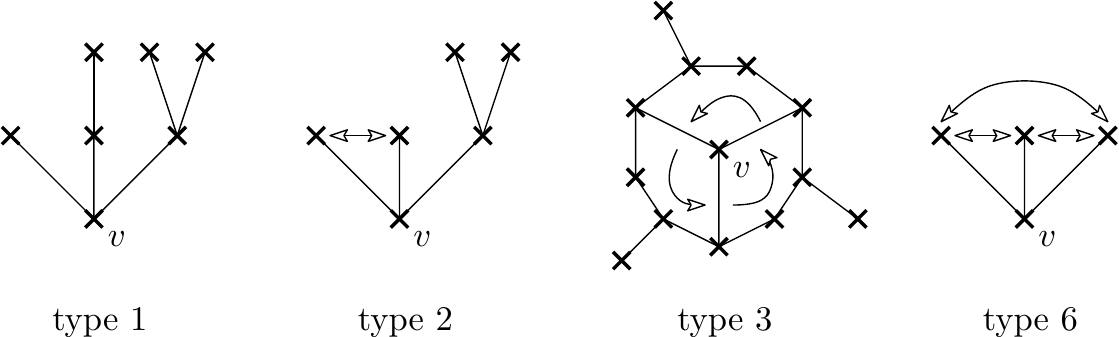}

	\caption{Examples of the different types of vertices.} \label{fig:vertices_type}
\end{figure}

The general plan is to understand distinguishability of cubic or subcubic graphs by looking at the way $\Aut_v(G)$ acts on the edges incident to a vertex $v$ of valence 3. If that action is trivial, call $v$ {\em type 1}. If the action leaves one edge fixed but interchanges the other two edges, call it {\em type 2}. If it fixes no edge, but is not  $S_3$, call it {\em type 3} and {\em type 6} otherwise. We note that $Aut_v(G)$ defines a permutation group $A$ on the neighbors of $v$ and the type is the same as the order of $A$.
See figure~\ref{fig:vertices_type} for some examples.

We can define similarly type 1 and type 2 for vertices of valence 2. All vertices of valence 1 are type 1.

\paragraph{Observation 1} If every valence 3 vertex of $G$ is type 3 or 6 and every valence 2 vertex is type 2, then $G$ is edge transitive.

\paragraph{Observation 2} In each of the five families $T_n, S_n, R_n^1,R_n^2,R_n^3$ with $D(G) = 3$, the root vertex is type 6 and all other vertices of valence 3 are type 2. For $S_n$,  the valence 2 vertices are type 1. For $R_n^1$, the valence 2 vertices are type 2.

\subsection{Vertices of type 1}

\begin{theorem} \label{thm:type1_val2}
	Suppose that $G$ has a valence 2 vertex $v$ of type 1. Then either $D(G) = 2$ or $G = S_n$ for some $n$.
\end{theorem}

\proof
Take a canonical coloring with $K$ the graph spanned by $v$ and its two neighbors $u$ and $w$.  By Proposition Internal (\ref{prop:internal}), any non-identity, color-preserving automorphism must move $v$ to another vertex $x$, which is either internal to  $K$ or a kiwi vertex of $S_1(K)$.

	If $x$ is internal to $K$, then $K$ must be a triangle with two vertices of valence 2. 
	Let $G'$ be obtained by removing all such triangles. If $D(G') = 2$, then we can extend any distinguishing coloring of $G$ to $G'$ by coloring the two valence 2 vertices in each such triangle, one white and one black. Therefore $D(G') = 3$, forcing $G' = T_n$, so $G$ is $S_n$.

Suppose instead that $x$ is a kiwi vertex in $S_1(K)$. Then by Proposition Internal (\ref{prop:internal}), the neighbors of $x$ are $u$ and $w$.  But $x$ would only be black if there is an automorphism fixing $u, w$ and interchanging $x$ with some other $z$. Then $z$ has valence 2 as well and is adjacent to $u, w$, forcing $G = K_{2,3}$. But then $v$ is not a type 1 vertex.
\qed

\begin{figure}[h]
	\centering

	\includegraphics[scale=0.8]{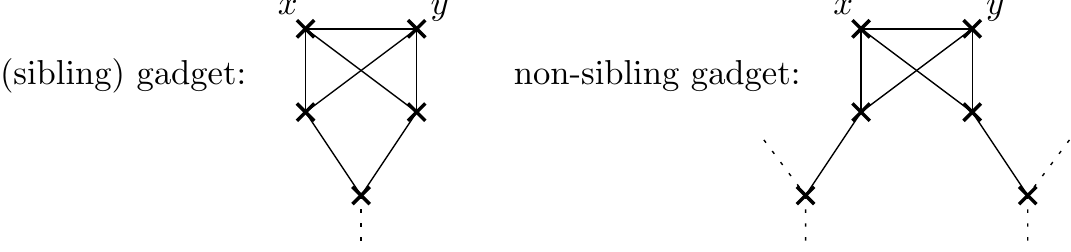}

	\caption{Example of a (sibling) gadget and a non-sibling gadget.} \label{fig:non_sibling_gadget}
\end{figure}

In every figure accompanying the definition of a gadget there are vertices $u$ and $v$. If $u$ and $v$ are not siblings, call the corresponding gadget a {\em non-sibling gadget}. See Figure \ref{fig:non_sibling_gadget} for an example.

\begin{corollary} \label{cor:nonsibling}
	If $G$ contains a non-sibling gadget, then $D(G) = 2$.
\end{corollary}

\proof
	Suppose first that $u,v$ are adjacent valence 2 vertices that are not siblings. The last vertex in the path of valence 2 vertices containing $u,v$ (in either direction) is a valence 2 vertex of type 1. Since $S_n$ does not contain a pair of adjacent vertices of valence 2 that are not siblings, $G \neq S_n$. Thus $D(G) = 2$ by Theorem~\ref{thm:type1_val2}.

	For each of the gadgets, replace all appearances of a non-sibling gadget by an edge to create a graph $G'$, with adjacent vertices of valence 2 that are not siblings. Then, by the above, $D(G') = 2$ and one easily extends a distinguishing 2-coloring from $G'$ to $G$, coloring $x,y$ black-white for $R_n^1, R_n^2$ and also $z,w$ black-white for $R_n^3$.
\qed

\begin {theorem} \label{thm:type1_val3}
	If $G$ has a type 1 vertex $v$ of valence 3, then $D(G) = 2$.
\end{theorem}

\proof
	Choose a canonical 2-coloring with $K$ spanned by $v$ and its three neighbors. This breaks all automorphisms in $\Aut_v(G)$. Thus the only color-preserving automorphisms left must move $v$ to a kiwi vertex $u$. By Proposition All Black (\ref{prop:all_black}), $u$ is internal to $K\cup S_1(K)$.

	Suppose that $u \neq v$ is internal to $K$. This forces $K$ to be a sibling or non-sibling type 2 gadget, making $D(G) = 2$ or $G = R_n^2$ by Corollary \ref{cor:nonsibling} and Lemma \ref{lem:sibling}. Suppose instead that $n = 1$ and $x_1$ is kiwi, which forces its down neighbors to be $u_0, w_0, z_0$. But $x_1$ would only be black if there was an automorphism fixing $K$ and interchanging $x_1$ with another $y_1$. This forces $G = K_{3,3}$, contradicting that $v$ has type 1.
\qed

\subsection{Type 2 vertices of valence 2}

\begin{theorem} \label{thm:type2_val2}
	If $G$ has a valence 2 vertex of type 2 but none of type 1, then  $G = K_{2,3}, R^1_n$ or $D(G) = 2$.
\end{theorem}

\proof
	Let $G'$ be the cubic graph obtained by smoothing over all valence 2 vertices. Thus $G$ is obtained from $G'$ by inserting valence 2 vertices in some edges. Note that we cannot have more than one such vertex in any edge of $G'$, since otherwise along this edge there will be a type 1 vertex of valence 2 in $G$.

	A multiple edge in $G'$, comes from a gadget of type 1 or a non-sibling gadget of type 1 or from $G = K_{2,3}$. Thus by Corollary~\ref{cor:nonsibling} either $D(G) = 2$ or $G = R_n^1$. We therefore assume that $G'$ has no multiple edges. If $D(G') = 2$, then $D(G) = 2$.
	Otherwise, either $G = R_n^2$ or $G = R_n^3$, since $G'$ is cubic.
	Since all edges of $G'$ except those in the gadgets have no automorphism interchanging the endpoints, the inserted vertices being type 2 must be in the gadget edges.

	Color one gadget vertex black and the rest of the valence 2 vertices white. In effect, this fixes one leaf $w$ of $T_n$. Now canonically color $T_n$ rooted at the center $v$ so that the neighbor of $v$ in the branch containing $w$ is colored white and the other two neighbors are colored black-white. This fixes the neighbors of $v$ and hence breaks all automorphisms. We conclude that $D(G') = 2$ and so $D(G) = 2$.
\qed

At this point, our classification is complete for subcubic graphs.

\subsection{Type 2 vertices of valence 3}

\begin{lemma} \label{lem:K_2,3}
	If the cubic graph $G$ contains $K_{2,3}$ as a subgraph, then
	$G = K_{3,3}, R_n^3$ or $D(G) = 2$.
\end{lemma}

\proof
	We note that the subgraphs of the gadget of type 3 spanned by $\{x,w ,u, y, z\}$ and $\{y,z,v,x,w\}$ are both isomorphic to $K_{2,3}$. We claim this the only way two copies of $K_{2,3}$ in $G$ can overlap. We view their union $H$ as obtained from two copies of $K_{2,3}$ with vertices identified in pairs. We consider  three different cases.

\begin{figure}[h]
		\centering
		
		\includegraphics[scale=0.9]{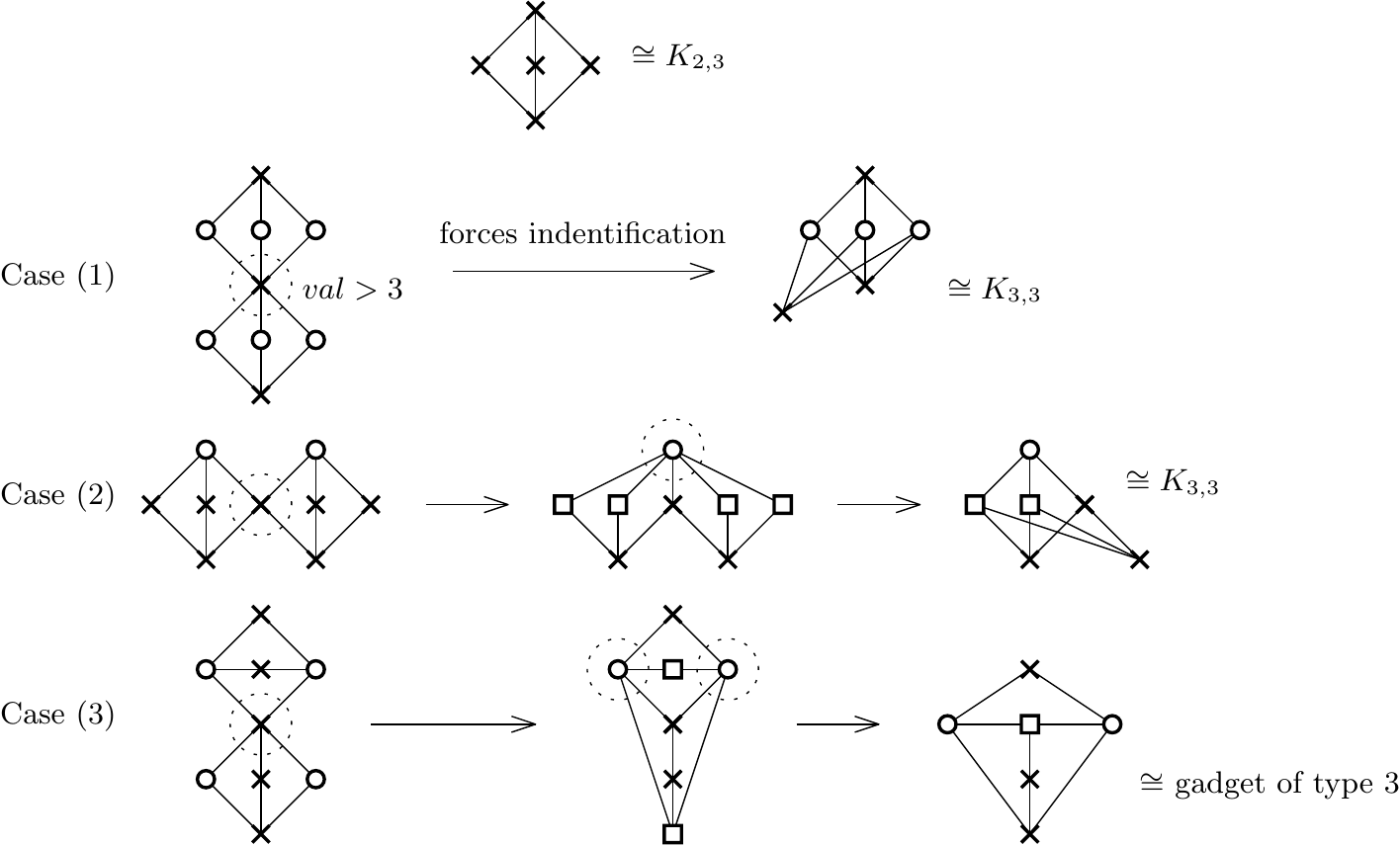}
		
		\caption{Identification of vertices.} \label{fig:identifications}
	\end{figure}

	\begin{itemize}
		\item \textit{Case 1:} Identifying a pair of valence 3 vertices also identifies in pairs their 3 neighbors (so the resulting vertex has valence 3); this yields $K_{3,3}$. See Figure~\ref{fig:identifications} for  illustration.
		\item \textit{Case 2:} Identifying two vertices of valence 2 forces the identification of a pair of their neighbors of valence 3, which was just  considered (see Case 1).
		\item \textit{Case 3:} Any identification of a valence 2 with a valence 3 vertex forces the identification of two other neighboring pairs of neighboring vertices, which in turn forces further identification. Thus $H$ has at most $6$ vertices and must be obtained by adding a single vertex of valence 2 to $K_{2,3}$, yielding gadget 3.
	\end{itemize}
	
	If $G$ contains a gadget of type 3, then $D(G) = 2$ or $G = R_n^3$. Therefore we assume that all copies in $G$ of $K_{2,3}$ are disjoint. Let $G'$ be the graph obtained by identifying in each copy of $K_{2,3}$ the two valence 3 vertices. The resulting graph has vertices of valence 2 with some vertices of valence 3 surrounded by vertices of valence  $2$.
	Thus $G'$ is not $S_n$ or $R_n^1$, so $D(G') = 2$. Any distinguishing coloring of $G'$ can be extended to one of $G$ by coloring the two valence 3 vertices of each $K_{2,3}$ black and white.
\qed

\begin{theorem} \label{thm:type2_val3}
	If the cubic graph $G$ has a type 2 vertex of valence 3, then $D(G) = 2$ or $G$ is $R_n^2$ or $R_n^3$. If $G$ has a type 3 vertex, then $D(G) = 2$.
\end{theorem}

\proof
	Let $v$ be a type 2 vertex of valence 3 and let $u_1, v_1$ be its neighbors which can be interchanged by $\Aut_v(G)$. Color $v$ black and all its neighbors black as in a canonical coloring rooted at the graph spanned by $v$ and its neighbors. Suppose that $u_1$ and $v_1$ are adjacent.  If their up neighbors are different, color one white and one black. This will also fix $u_1, v_1$ and can be continued to a canonical coloring that breaks all automorphisms in $\Aut_v(G)$. If instead $u_1, v_1$ have a common up neighbor, we have a gadget 2 or a non-sibling gadget 2.  Then either $D(G) = 2$ or $G = R_n^2$ by Lemma~\ref{lem:sibling} and Corollary \ref{cor:nonsibling}.

	We assume therefore that $u_1$ and $v_1$ are not adjacent. If they have one common up neighbor but not two, color that neighbor white and the other two up neighbors black and white. This distinguishes $u_1, v_1$. Suppose $u_1, v_1$ have two common up neighbors. Then $G$ contains a $K_{2,3}$ so by Lemma~\ref{lem:K_2,3} we infer that  $D(G) = 2$ or $G = R^3_n$.

\begin{figure}[h]
		\centering
		
		\includegraphics[scale=0.9]{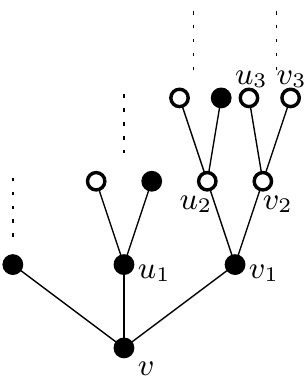}
		
		\caption{Example of the coloring to fix the vertices $u_1$ and $v_1$ by any automorphism that fixes $v$.} \label{fig:thm_type2_val3}
	\end{figure}

	Therefore, assume that the up neighbors of $u_1$ and $v_1$ are distinct. Color the up neighbors of $u_1$ black-white and the up neighbors $u_2, v_2$ of $v_1$ both white. This distinguishes $u_1, v_1$ but allows an interchange of $u_2, v_2$.
	Repeat this process on $u_2, v_2$. Either we distinguish $u_2, v_2$ or we find $u_3, v_3$ that can be interchanged. Continue the process. If $G$ is finite, the process must end either with
	$G = R_n^2, R_n^3$
	or with a 2-coloring which breaks all automorphisms in $\Aut_v(G)$ with $v$ and its three neighbors colored black. If $G$ is infinite, we continue the process as long as $u_n, v_n$ have distinct up neighbors, yielding a 2-coloring that breaks any automorphism fixing $v$ see Figure~\ref{fig:thm_type2_val3}.

	We then can proceed as in the proof of Proposition All Black (\ref{prop:all_black}) to show that there is no other black vertex surrounded by three black vertices for $S_n(K)$, $n > 1$, and by construction there is none in $S_1(K)$. We conclude that no color-preserving automorphism moves $v$, so the coloring is distinguishing.

	Suppose instead that $v$ is a type 3 vertex. We proceed exactly as before except now we are breaking any automorphism taking $u_1$ to $v_1$ or $v_1$ to $u_1$ (but not interchanging them). We still break all non-identity elements of $\Aut_v(G)$. Since neither $R_n^2$ or $R_n^3$ have a type 3 vertex, we must have $D(G) = 2$.
\qed

We have now completed the classification for graphs $G$ with a vertex of type 1,2,3. There remains only the case where all vertices have type 6. Then $G$ must be cubic and, as is easily seen, edge transitive. We will treat the distinguishability of edge transitive graphs in the next section.

\section {Girth}

Our analysis of edge transitive cubic graphs $G$ uses the girth of $G$.
One easily verifies that the only edge transitive cubic graph of girth 3 is $K_4$, and that $K_{3,3}$ and the cube are the only edge transitive graphs of girth 4. We know $D(K_4) = D(K_{3,3}) = 4$ and $D(Q) = 3$.

For girth 5 we observe that edge transitive graphs that are not vertex transitive must be bipartite.
Hence all edge transitive graphs of odd girth are also vertex transitive. But there are only two vertex transitive cubic graphs of girth 5, the dodecahedron $H$ and the Petersen graph \cite{GM07}. In Lemma \ref{lem:P} we have shown that $D(P) = 3$. However, $D(H) = 2$. To find a 2-distinguishing coloring color black a vertex $v$, its three neighbors $x,y,z$, and a vertex $w$ adjacent to $x$.

For girth $s > 5$ the situation changes drastically. Although there are  only five cubic edge transitive graphs of girth at most 5, and no infinite ones, there are infinitely many such graphs already for girth 6, and an infinite one is the honeycomb lattice, which is also edge transitive.

Thus it remains for us to show that edge transitive cubic graphs of girth $s > 5$ are 2-distinguishable. In fact, our proof does not use edge transitivity at all.

\begin{lemma} \label{lem:girth_over_7}
	If $G$ is a cubic graph with girth $s > 6$, then $D(G) = 2$.
\end{lemma}

\begin{figure}[h]
		\centering
		
		\includegraphics[scale=0.9]{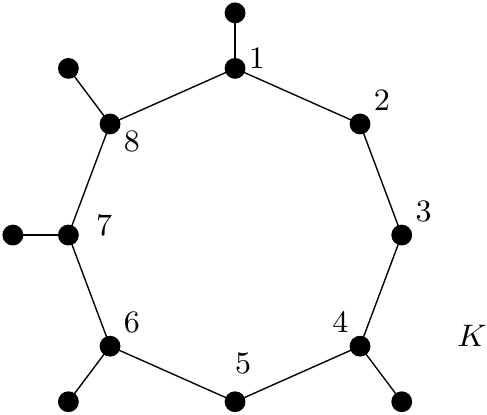}
		
		\caption{The set $K$ for $s=8$. The vertices $1$ and $4$ together with $6$, $7$ and $8$ are the kiwi vertices.} \label{fig:circle_with_whiskers}
	\end{figure}

\proof
	Let $C$ be a cycle of length $s$. Since $s > 6$, 
	each vertex in $S_1(C)$ is adjacent to only one vertex in $C$. Moreover, if two vertices in $S_1(C)$ are adjacent, then they can be used to form a path of length three between two vertices in $C$ of distance at most $s/2$, contradicting $s > 6$. Let the vertices of $C$ be denoted by $1, 2, 3, \dots, s$. Let $K$ be $C$ together with the whiskers at vertices $1$ and $4$ as well as $6, \dots, s$, see Figure~\ref{fig:circle_with_whiskers}.

	Choose a canonical 2-coloring rooted at $K$. By Proposition All Black (\ref{prop:all_black}), there is no kiwi vertex in  $S_n(K)$
	for $n > 1$. There is none in $S_1(K)$ either, since any such vertex would be adjacent to three vertices in $K$, forcing a $K_{2,3}$ in $G$, contradicting $s > 6$.
	Thus any color-preserving automorphism $\phi$ must leave invariant the kiwi vertices $1, 4, 6, \dots, s$. The graph spanned by these vertices consists of an isolated vertex $4$ and a path $6, \dots, s, 1$. Thus $\phi$ fixes $4$, and either leaves the path fixed or reverses it (interchanging vertices $1,6$). In the first case, $\phi$ fixes $5$, since $s > 4$, and $\phi$ fixes $2,3$ since $1,\phi(2), \phi(3), 4$ provides another path of length 3 between 1 and 4, contradicting $s > 6$. In the second case, $4, \phi(5), \phi(6) = 1$ provides a path of length $2$ from $4$ to $1$, contradicting $s > 5$. We conclude that $\phi$ fixes all vertices of $C$ and hence all vertices of $K$, so $\phi$ is the identity.
\qed

For girth $s = 6$, we also have $D(G) = 2$ but the argument is slightly more complicated.

\begin{lemma} \label{lem:girth6}
	Let $G$ be a cubic graph with girth $s = 6$. Then $D(G) = 2$.
\end{lemma}

\proof
	We note that the proof of Lemma~\ref{lem:girth_over_7} for girth $s > 6$ only uses $s \neq 6$ to insure that for a cycle $C$ of length $s$, the graph $K$ obtained from $C$ with whiskers to vertices $1', 4', 6', ..., s'$ has no edges between the whiskers. The rest of the proof only requires $s > 5$. In particular, for girth $s = 6$, if $G$ contains a cycle $123456$ with no edge between $1'$ and $4'$, then $D(G) = 2$ (as there can be no edge between $4'$ and $6'$ as otherwise $s < 6$). By cyclically permuting $1,2,3,4,5,6$, we conclude that for every 6-cycle $C$, we must have edges $1'4', 2'5', 3'6'$ in $S_1(C)$ or else $D(G) = 2$. Applying this to the cycle $1'12344'$, we have that $1'$ and $3'$ must have a common neighbor. Since the choice of which vertex is labeled $1$ does not matter, $3'$ and $5'$ have a common neighbor, as do $5'$ and $1'$.  Since all vertices have valence 3, it must be that $1',3',5'$ have one common neighbor $7$. Similarly, $2',4',6'$ have one common neighbor $8$.

\begin{figure}[h]
		\centering
		
		\includegraphics[scale=0.9]{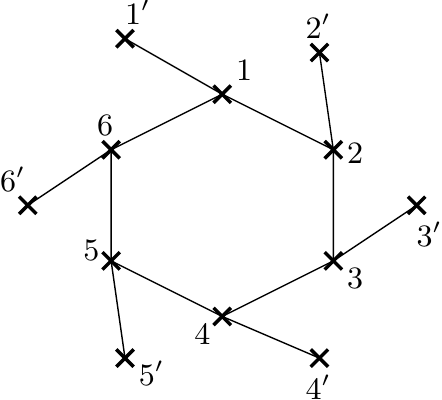}
		\includegraphics[scale=0.9]{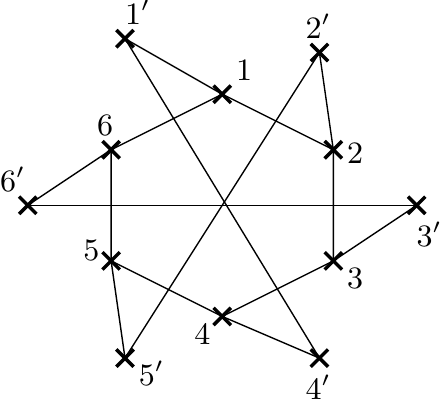}
		\includegraphics[scale=0.9]{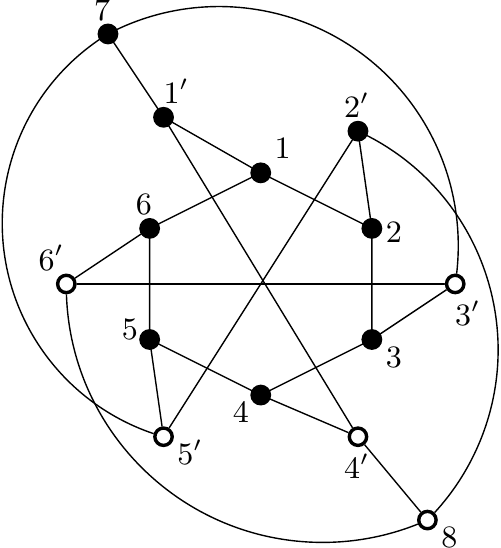}
		
		\caption{Construction of the graph described  in the proof of Lemma~\ref{lem:girth6}.} \label{fig: girth6}
	\end{figure}

	At this point all 14 of the vertices have valence 3 so we have the entire graph, see Figure \ref{fig: girth6}. This is the Heawood graph (the dual of the triangulation of the torus with underlying graph $K_7$).

	Now consider the following 2-coloring of the graph. Let $1, 2, 3, 4, 5, 6, 1', 2', 7$ be black and the remaining vertices white.
	In the graph $H$ spanned by the black vertices, $7$ is the only vertex of valence 1 adjacent to a vertex of valence 2 (namely $1'$). Thus any color-preserving automorphism $\phi$ fixes $7$. Thus $\phi$ also fixes $1'$ and hence $1$. Since $2'$ is the only vertex of valence one, $\phi$ also fixes $2$. Thus $\phi$ fixes the remaining vertices of the cycle $C$, so $\phi$ fixes all black vertices. But then $\phi$ also fixes the white vertices adjacent to $3, 4, 5, 6$. That leaves only $8$ so it must be fixed as well, making $\phi$ the identity
	(compare the right graph in Figure \ref{fig: girth6} for the coloring).
\qed

\section {Questions}

There are a  number of questions worth further study,

\paragraph{Question 1} (Higher Valence) Can we classify graphs $G$ with $\Delta(G) = d = D(G)$?

As we have observed, if $G = T(n,d)$, then $\Delta(G) = d = D(G)$. We could add edges within each sibling family of size $d - 1$ to form a graph $S(n,d)$ analogous to $S_n$ (the vertices of a sibling family then have valence $d-1$). We can also attach $d - 1$ independent vertices to a sibling family using $K_{d-1, d-1}$ to obtain a graph analogous to $R_n^1$. There do not appear to be analogues for $R_n^2$ and $R_n^3$.

We can define a canonical $d - 1$ coloring rooted at a graph $K$ such that the only color-preserving automorphism of $G$ fixing the vertices of $K$ is the identity. Then we have to identify properties of such colorings that restrict the structure of $K$ and $S_1(K)$. Note that a variation of the canonical coloring using $d+1$ colors, with color $d+1$ for a vertex $v$ , colors $1, 2, \dots, d$ for the neighbors of $v$, and colors $1, 2, \dots, d-1$ is how one gets $D(G) \leq d+1$. And to show $D(G)=d+1$ only for $G=K_{d+1}$ or $G=K_{d,d}$, one uses canonical $d$-colorings rooted at an asymmetric vertex-induced subgraph $K$ with color $d$ used only on $K$.

\paragraph{Question 2} (Highly Symmetric Graphs) If $\Delta(G) = d$ and $G$ is vertex transitive, must $D(G) = 2$ for all but finitely many $G$?

\paragraph{Question 3} (Connectivity) What is the relationship between vertex or edge connectivity, valence, and distinguishing number?

The examples with $D(G) = d$ are not 2-connected. What happens if we require, say, 3-connectivity? For example, we can get $D(G) = d - 3$ with $G$ 3-connected by attaching $K_{d-1}$ at three vertices of valence 2. As the connectivity goes up, the distinguishing number seems to go down, with finitely many exceptions like $K_{d+1}$.

\paragraph{Question 4} (Infinite Graphs) What happens for infinite $G$ with $\Delta(G) = d > 3$?

For infinite graphs, we expect that if $\Delta(G) = d$, then $D(G) < d$, just as for $d = 3$. But there are also interesting questions just for such $G$ with $D(G)=2$. As we observed before Corollary 2.6, for finite graphs, large enough motion implies $D(G)=2$.  The Infinite Motion Conjecture \cite{T11} 
is that if $G$ is locally finite and $m(G) = \infty$, then $D(G) = 2$.  On the other hand, for the case $d = 3$, we see there is no need for the hypothesis of infinite motion to get $D(G) = 2$, and there are other classes of graphs with $D(G) = 2$ that do not depend intrinsically on infinite motion \cite{STW12}. As we observed, however, it is easy to construct an infinite $d$-valent graph $G$ with $D(G) = d-1$, so for $d > 2$, we expect infinite motion to be involved.

\paragraph{Question 5} (Motion).  For cubic graphs, if the motion $m(G)>2$, then $D(G)=2$ with the exception of $Q$ and $P$.  For $d>3$, is it the case that if $m(G)>d$, then $D(G)=2$ with finitely many exceptions?  

Perhaps, even $m(G)>2$ gives $D(G)=2$ with finitely many exceptions.  

\paragraph{Question 6} (Chromatic Distinguishing Number) Suppose all colorings are required to be proper (adjacent vertices get different colors).  What happens when $\Delta(G) = 3$?

Collins and Trenk \cite{CT06} define the chromatics distinguishing number $\chi_D(G)$ to be the least $k$ such that $G$ has a proper $k$-coloring whose only color-preserving automorphism is the identity. They prove that $\chi_D(G) \leq 2d$ with equality only for $K_{d,d}$ and $C_6$. For $d = 3$, there is the possibility of classifying graphs with $D(G) = 5$, especially in the case that $G$ is bipartite.

In \cite{ikps-17} the chromatic distinguishing number of infinite graphs is investigated. For connected graphs of bounded valence $d$ it is shown that $\chi_D(G) \leq 2d-1$, and for infinite subcubic graphs of  infinite motion this improves to $\chi_D(G) \leq 4$.

\paragraph{Question 7} (Edge Distinguishing)
One can also define \cite{KP15} the \emph{distinguishing index} (or edge distinguishing number)  $D'(G)$ as the least $k$ such that some $k$-coloring of the edges of $G$ is preserved only by the identity. In \cite{KP15} it is  shown that $D'(G) \leq \Delta(G)$ for finite graphs. For infinite graphs $\Delta(G)$ has to be replaced by the supremum of the valences \cite{BP15}.

What happens with $D'(G)$ when $\Delta(G) = 3$?

\paragraph{Question 8} (Cost)  When $D(G) = 2$, the cost \cite{B13,BI17} is the least number of times the color black is used. When $\Delta(G) = D(G) = 2$, what can we say about the cost? For  cubic graphs this is treated in \cite{ikps-17}.

The canonical coloring tends to use black as few times as possible for $S_n(K)$, $n > 0$.  How close does this number come to the cost?

\end{document}